\documentclass[12pt]{article}%
\usepackage{amsmath}
\usepackage{amsfonts}
\usepackage{amssymb}
\usepackage{graphicx}%
\providecommand{\U}[1]{\protect\rule{.1in}{.1in}}
\newtheorem{theorem}{Theorem}
\newtheorem{acknowledgement}[theorem]{Acknowledgement}

\begin{document}

\title{Asymptotic analysis of a family of polynomials associated with the inverse
error function}
\author{Diego Dominici\\Department of Mathematics\\State University of New York at New Paltz\\1 Hawk Dr. Suite 9\\New Paltz, NY 12561-2443\\dominicd@newpaltz.edu
\and Charles Knessl\\Department of Mathematics, Statistics and Computer Science \\University of Illinois at Chicago (M/C 249) \\851 South Morgan Street \\Chicago, IL 60607-7045\\knessl@uic.edu}

\maketitle

\begin{abstract}
We analyze the sequence of polynomials defined by the differential-difference equation $P_{n+1}(x)=P_{n}^{\prime}(x)+x(n+1)P_{n}(x)$ asymptotically as $n\rightarrow\infty$. The polynomials $P_{n}(x)$ arise in the computation of higher derivatives of the inverse error function $\operatorname{inverf}(x)$. We use singularity analysis and discrete versions of the WKB and ray methods and give numerical results showing the accuracy of our formulas.
\end{abstract}

MSC-Class: 33B20 (Primary) 34E20, 33E30 (Secondary).


\section{Introduction}

The error function $\operatorname{erf}(x)$ is defined by \cite{MR1225604}
\begin{equation}
\operatorname{erf}(x)=\frac{2}{\sqrt{\pi}}\int\limits_{0}^{x}\exp\left(
-t^{2}\right)  dt \label{erf}%
\end{equation}
and its inverse $\operatorname*{inverf}\left(  x\right)  $, which we will
denote by $\mathfrak{I}(x),$ satisfies $\mathfrak{I}\left[  \operatorname{erf}%
(x)\right]  =\operatorname{erf}\left[  \mathfrak{I}(x)\right]  =x.$ The
function $\mathfrak{I}(x)$ appears in several problems of heat conduction
\cite{MR0434172}. In \cite{MR1986919} we considered the function
\begin{equation}
N(x)=\frac{1}{\sqrt{2\pi}}\int_{-\infty}^{x}e^{-t^{2}/2}dt \label{N}%
\end{equation}
and its inverse $S(x),$ satisfying
\[
S\left[  N(x)\right]  =N\left[  S(x)\right]  =x.
\]
It is clear from (\ref{erf}) and (\ref{N}) that%
\[
N(x)=\frac{1}{2}\left[  \operatorname{erf}\left(  \frac{x}{\sqrt{2}}\right)
+1\right]
\]
and therefore%
\begin{equation}
S(x)=\sqrt{2}\mathfrak{I}\left(  2x-1\right)  . \label{SJ}%
\end{equation}

In \cite{MR1986919} we showed that%
\begin{equation}
S^{\prime}(x)=\sqrt{2\pi}\exp\left[  \frac{1}{2}S^{2}(x)\right]  \label{ds}%
\end{equation}
and
\begin{equation}
S^{(n)}=P_{n-1}(S)(S^{\prime})^{n}\quad n\geq1, \label{Sn}%
\end{equation}
where $P_{n}(x)$ is a polynomial of degree $n$ satisfying the recurrence
\begin{equation}
P_{0}(x)=1,\quad P_{n+1}(x)=P_{n}^{\prime}(x)+x(n+1)P_{n}(x),\quad n\geq1.
\label{2.3}%
\end{equation}
The same approach was employed by Carlitz in \cite{MR0153878}. From
(\ref{2.3}), it follows easily that for a fixed value of $n$
\begin{equation}
P_{n}(x)\sim n!\,x^{n},\quad x\rightarrow\infty. \label{large x}%
\end{equation}

From (\ref{SJ}) and (\ref{Sn}), we conclude that%
\[
\mathfrak{I}^{(n)}=2^{\frac{n-1}{2}}P_{n-1}\left(  \sqrt{2}\mathfrak{I}%
\right)  \left(  \mathfrak{I}^{\prime}\right)  ^{n}\quad n\geq1.
\]
Since%
\[
\mathfrak{I}(0)=0,\quad\mathfrak{I}^{\prime}(0)=\frac{\sqrt{\pi}}{2},
\]
we have%
\begin{equation}
\mathfrak{I}^{(n)}(0)=\frac{1}{\sqrt{2}}\left(  \frac{\pi}{2}\right)
^{\frac{n}{2}}P_{n-1}(0). \label{Ider}%
\end{equation}
It follows from (\ref{Ider}) that estimating $\mathfrak{I}^{(n)}(0)$ for large
values of $n$ is equivalent to finding an asymptotic approximation of the
polynomials $P_{n}(x)$ when $x=0$.

The objective of this work is to study $P_{n}(x)$ asymptotically as
$n\rightarrow\infty$ for various ranges of $x.$ We shall obtain different
asymptotic expansions for $n\rightarrow\infty$ and (i) $0<x<\infty,$ (ii)
$x=O\left(  n^{-1}\right)  $ and (iii) $x=O\left(  \sqrt{\ln(n)}\right)  .$
The paper is organized as follows: in Section \ref{section2} we approach the
problem using a singularity analysis of the generating function
\cite{MR2078410} of the polynomials $P_{n}(x).$ In Section \ref{section3} we
apply the WKB method to the differential-difference equation (\ref{2.3}). In
\cite{MR1695194}, we used this approach in
the asymptotic analysis of computer science problems and in \cite{MR2390273}
to study the Krawtchouk polynomials. Finally, in Section \ref{section4} we
analyze (\ref{2.3}) again using the ray method \cite{MR499726} and obtain an
asymptotic approximation valid in various regions of the $(x,n)$ domain. In \cite{MR2364955},
\cite{MR2401156}, \cite{MR2250392}, we employed the same
technique to analyze asymptotically other families of polynomials and in \cite{MR2117327}, 
\cite{MR2262956} to study some queueing problems.

\section{Singularity analysis\label{section2}}

In \cite{MR1986919} we obtained the exponential generating function
\[
\sum_{n=0}^{\infty}P_{n}(x)\frac{z^{n}}{n!}=\exp\left\{  \frac{1}{2}%
S^{2}[N(x)+zN^{\prime}(x)]-\frac{x^{2}}{2}\right\}  ,
\]
which implies that%
\[
P_{n}(x)=e^{-x^{2}/2}\frac{n!}{2\pi\mathrm{i}}\oint\limits_{\left\vert
z\right\vert <r}\exp\left\{  \frac{1}{2}S^{2}[N(x)+zN^{\prime}(x)]\right\}
\frac{dz}{z^{n+1}},
\]
where the integration contour is a small loop around the origin in the complex
plane. Using (\ref{ds}), we have%
\begin{align*}
P_{n}(x)  &  =e^{-x^{2}/2}\frac{n!}{2\pi\mathrm{i}}\oint\limits_{\left\vert
z\right\vert <r}\frac{1}{\sqrt{2\pi}}S^{\prime}[N(x)+zN^{\prime}(x)]\frac
{dz}{z^{n+1}}\\
&  =e^{-x^{2}/2}\frac{1}{\sqrt{2\pi}N^{\prime}(x)}\frac{n!}{2\pi\mathrm{i}%
}\oint\limits_{\left\vert z\right\vert <r}\frac{1}{z^{n+1}}dS[N(x)+zN^{\prime
}(x)]\\
&  =\frac{n!}{2\pi\mathrm{i}}\oint\limits_{\left\vert z\right\vert <r}%
\frac{n+1}{z^{n+2}}S[N(x)+zN^{\prime}(x)]dz
\end{align*}
and therefore%
\begin{equation}
P_{n}(x)=\frac{\left(  n+1\right)  !}{2\pi\mathrm{i}}\oint\limits_{\left\vert
z\right\vert <r}S[N(x)+zN^{\prime}(x)]\frac{dz}{z^{n+2}}. \label{Pn1}%
\end{equation}

Since $S(x)$ has singularities at $x=0$ and $x=1,$ we consider the functions%
\[
Z_{1}(x)=\frac{1-N(x)}{N^{\prime}(x)},\quad Z_{0}(x)=-\frac{N(x)}{N^{\prime
}(x)}.
\]
We have%
\begin{equation}
Z_{1}(-x)=-\sqrt{2\pi}\exp\left(  \frac{x^{2}}{2}\right)  -Z_{1}(x)
\label{Z1-}%
\end{equation}
and%
\begin{align*}
Z_{1}(x)  &  =x^{-1}+O\left(  x^{-3}\right)  ,\quad x\rightarrow\infty,\\
\quad Z_{0}(x)  &  =-\sqrt{2\pi}\exp\left(  \frac{x^{2}}{2}\right)
+x^{-1}+O\left(  x^{-3}\right)  ,\quad x\rightarrow\infty.
\end{align*}
Changing variables to
\[
w=\left[  z-Z_{1}(x)\right]  N^{\prime}(x)
\]
in (\ref{Pn1}), we obtain%
\[
P_{n}(x)=\frac{1}{N^{\prime}(x)}\frac{\left(  n+1\right)  !}{2\pi\mathrm{i}%
}\oint\limits_{C}S\left(  w+1\right)  \left[  \frac{w}{N^{\prime}(x)}%
+Z_{1}(x)\right]  ^{-(n+2)}dw,
\]
or,
\begin{equation}
P_{n}(x)=\sqrt{\pi}e^{x^{2}/2}\frac{\left(  n+1\right)  !}{2\pi\mathrm{i}%
}\oint\limits_{C}\frac{\mathcal{J}\left(  w+1\right)  }{\left[  \sqrt
{\frac{\pi}{2}}e^{x^{2}/2}w+Z_{1}(x)\right]  ^{n+2}}dw, \label{contour}%
\end{equation}
where $C$ is a small loop about $w=w^{\ast}(x)$ in the complex plane, with
\[
w^{\ast}(x)=-\sqrt{\frac{2}{\pi}}e^{-x^{2}/2}Z_{1}(x).
\]

To expand (\ref{contour}) for $n\rightarrow\infty$ with a fixed $x\in\left(
0,\infty\right)  ,$ we employ singularity analysis. The function
$\mathcal{J(}w)$ has singularities at $w=\pm1.$ By (\ref{erf}), we have%
\[
w=\frac{2}{\sqrt{\pi}}\int\limits_{0}^{\mathcal{J}}e^{-t^{2}}%
dt=1-e^{-\mathcal{J}^{2}}\left[  \frac{1}{\sqrt{\pi}\mathcal{J}}+O\left(
\mathcal{J}^{-3}\right)  \right]  ,\quad\mathcal{J}\rightarrow\infty,
\]
so that%
\[
\mathcal{J(}w)\sim\sqrt{-\ln(1-w)},\quad w\rightarrow1^{-}%
\]
and by symmetry we have%
\[
\mathcal{J(}w)\sim-\sqrt{-\ln(1+w)},\quad w\rightarrow-1^{+}.
\]

\begin{figure}[t]
\begin{center}
\rotatebox{270} {\resizebox{10cm}{!}{\includegraphics{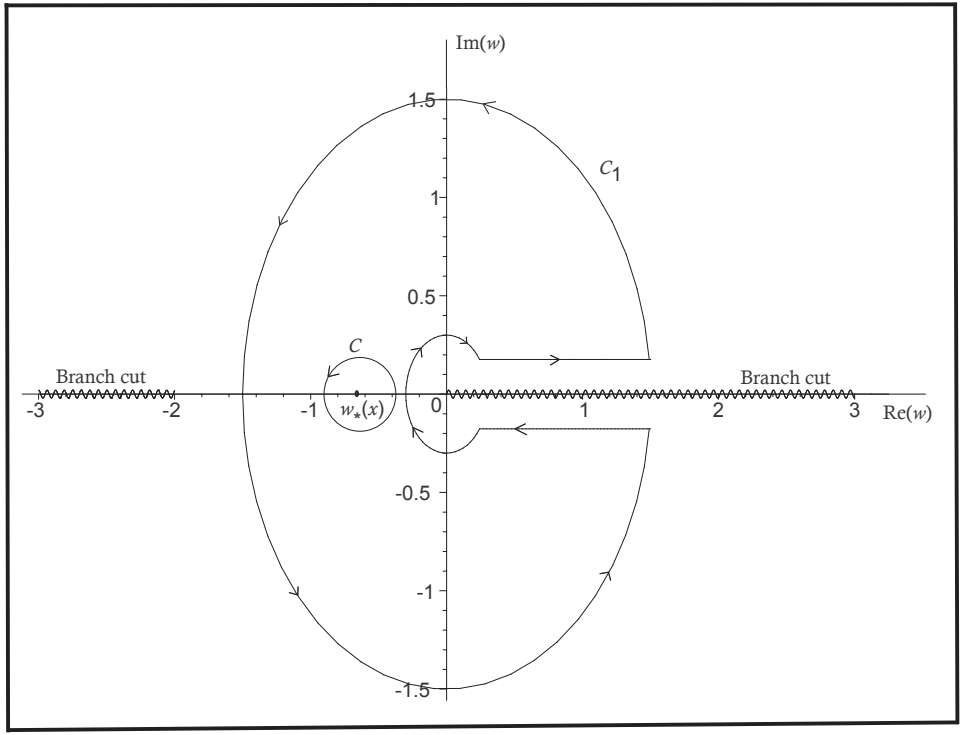}}}
\end{center}
\caption{A sketch of the contours $C$ and $C_{1}$.}%
\label{Fig12}%
\end{figure}

The integrand in (\ref{contour}) thus has singularities at $w=0$ and $w=-2,$
but for $x>0,$ the former is closer to $w^{\ast}(x).$ We expand (\ref{contour}%
) around $w=0$ by setting $w=\delta/n$ and using%
\begin{equation}
\left[  Z_{1}(x)+\sqrt{\frac{\pi}{2}}e^{x^{2}/2}\frac{\delta}{n}\right]
^{-(n+2)}\sim\left[  Z_{1}(x)\right]  ^{-(n+2)}\exp\left[  -\sqrt{\frac{\pi
}{2}}e^{x^{2}/2}\frac{\delta}{Z_{1}(x)}\right]  . \label{13}%
\end{equation}
Then, we deform the contour $C$ in (\ref{contour}) to a new contour $C_{1}$
that encircles the branch point at $w=0$ (see Figure \ref{Fig12}). This leads to%
\begin{align}
P_{n}(x)  &  \sim\left(  n+1\right)  !\sqrt{\pi}e^{x^{2}/2}\left[
Z_{1}(x)\right]  ^{-(n+2)}\frac{1}{n}\label{14}\\
&  \times\frac{1}{2\pi\mathrm{i}}\int\limits_{0}^{\infty}\left(  \Upsilon
^{+}-\Upsilon^{-}\right)  \exp\left[  -\sqrt{\frac{\pi}{2}}\frac{e^{x^{2}/2}%
}{Z_{1}(x)}\delta\right]  d\delta,\nonumber
\end{align}
where%
\[
\Upsilon^{\pm}(\delta,n)=\sqrt{\pm\mathrm{i}\pi-\ln\left(  \delta\right)
+\ln\left(  n\right)  }.
\]
Here $\Upsilon^{\pm}(\delta,n)$ corresponds to the approximation of
$\mathcal{J(}w+1)$ for $w\rightarrow0,$ above or below the right branch cut in
Figure \ref{Fig12}.

\begin{figure}[t]
\begin{center}
\rotatebox{270} {\resizebox{10cm}{!}{\includegraphics{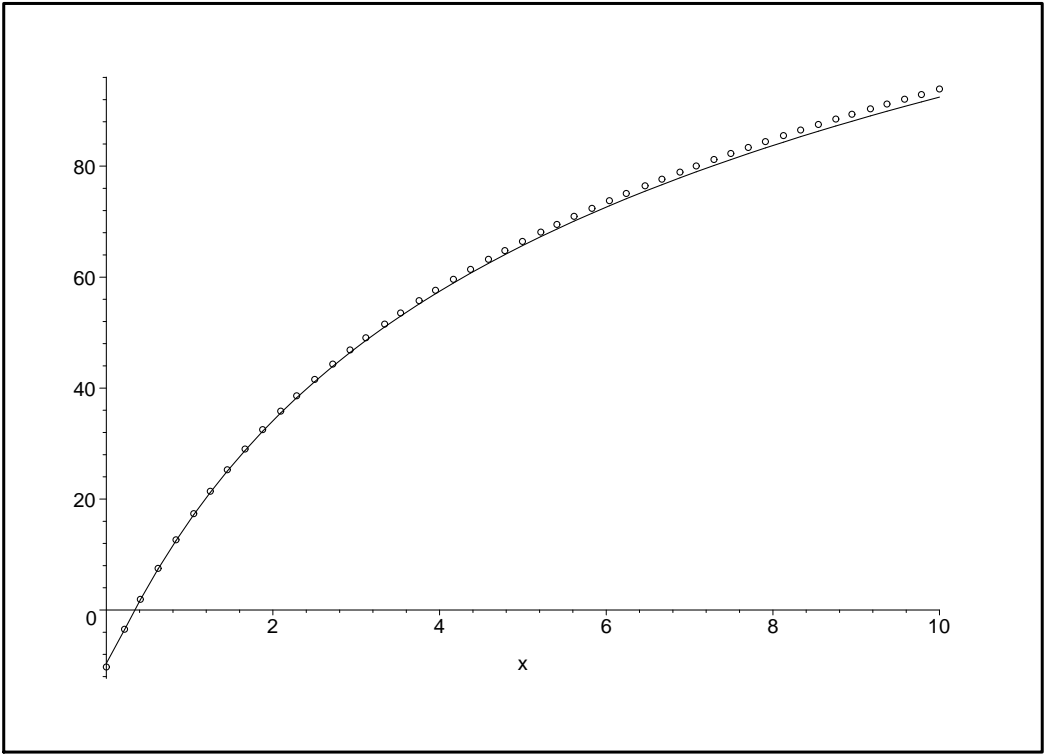}}}
\end{center}
\caption{A plot of $\ln\left[  P_{40}(x)/40!\right]  $ (solid line) and $\ln\left[  \Psi
_{1}(x,40)/40!\right]  $ (ooo).}%
\label{xO1}%
\end{figure}

For $n$ large we have%
\[
\Upsilon^{+}(\delta,n)-\Upsilon^{-}(\delta,n)\sim\frac{\pi\mathrm{i}}%
{\sqrt{\ln\left(  n\right)  }}%
\]
and then evaluating the elementary integral in (\ref{14}) leads to
$P_{n}(x)\sim\Psi_{1}(x,n)$ as $n\rightarrow\infty$ with%
\begin{equation}
\Psi_{1}(x,n)=\frac{n!}{\sqrt{2\ln\left(  n\right)  }}\left[  Z_{1}(x)\right]
^{-\left(  n+1\right)  }=\frac{n!}{\sqrt{2\ln\left(  n\right)  }}\left[
\frac{e^{-x^{2}/2}}{\zeta\left(  x\right)  }\right]  ^{n+1}\label{15}%
\end{equation}
and
\begin{equation}
\zeta\left(  x\right)  =\sqrt{\frac{\pi}{2}}\left[  1-\operatorname*{erf}%
\left(  \frac{x}{\sqrt{2}}\right)  \right]  \sim\exp\left(  -\frac{x^{2}}%
{2}\right)  \left[  x^{-1}+O\left(  x^{-3}\right)  \right]  ,\quad
x\rightarrow\infty.\label{zeta}%
\end{equation}
In Figure \ref{xO1} we plot $\ln\left[  P_{40}(x)/40!\right]  $ and $\ln\left[  \Psi
_{1}(x,40)/40!\right]  $. We see that the approximation is very good for
$x=O(1)$ but it becomes less precise as $x\rightarrow\infty.$ This is because
our previous analysis assumes that $n\rightarrow\infty$ with $0<x<\infty.$ If
either $x\rightarrow0$ or $x\rightarrow\infty,$ we must modify it, which we will do next.

When $x\rightarrow0,$ or more generally when $x=O\left(  n^{-1}\right)  ,$ the
singularities at $w=0$ and $w=-2$ are nearly equidistant from $w^{\ast}(x).$
On the scale $x=y/n,$ $y=O\left(  1\right)  ,$ we have%
\begin{equation}
Z_{1}(x)=Z_{1}\left(  \frac{y}{n}\right)  =\sqrt{\frac{\pi}{2}}-\frac{y}%
{n}+O\left(  n^{-2}\right)  ,\quad n\rightarrow\infty\label{16}%
\end{equation}
and (\ref{15}) simplifies to
\begin{equation}
P_{n}(x)\sim\frac{n!}{\sqrt{2\ln\left(  n\right)  }}\left(  \sqrt{\frac{2}%
{\pi}}\right)  ^{n+1}\exp\left(  y\sqrt{\frac{2}{\pi}}\right)  . \label{17}%
\end{equation}
But, to this we must add the contribution from $w=-2,$ which corresponds to
replacing $y$ by $-y$ and multiplying by $\left(  -1\right)  ^{n}$ the right
hand side of (\ref{17}).

\begin{figure}[t]
\begin{center}
\rotatebox{270} {\resizebox{10cm}{!}{\includegraphics{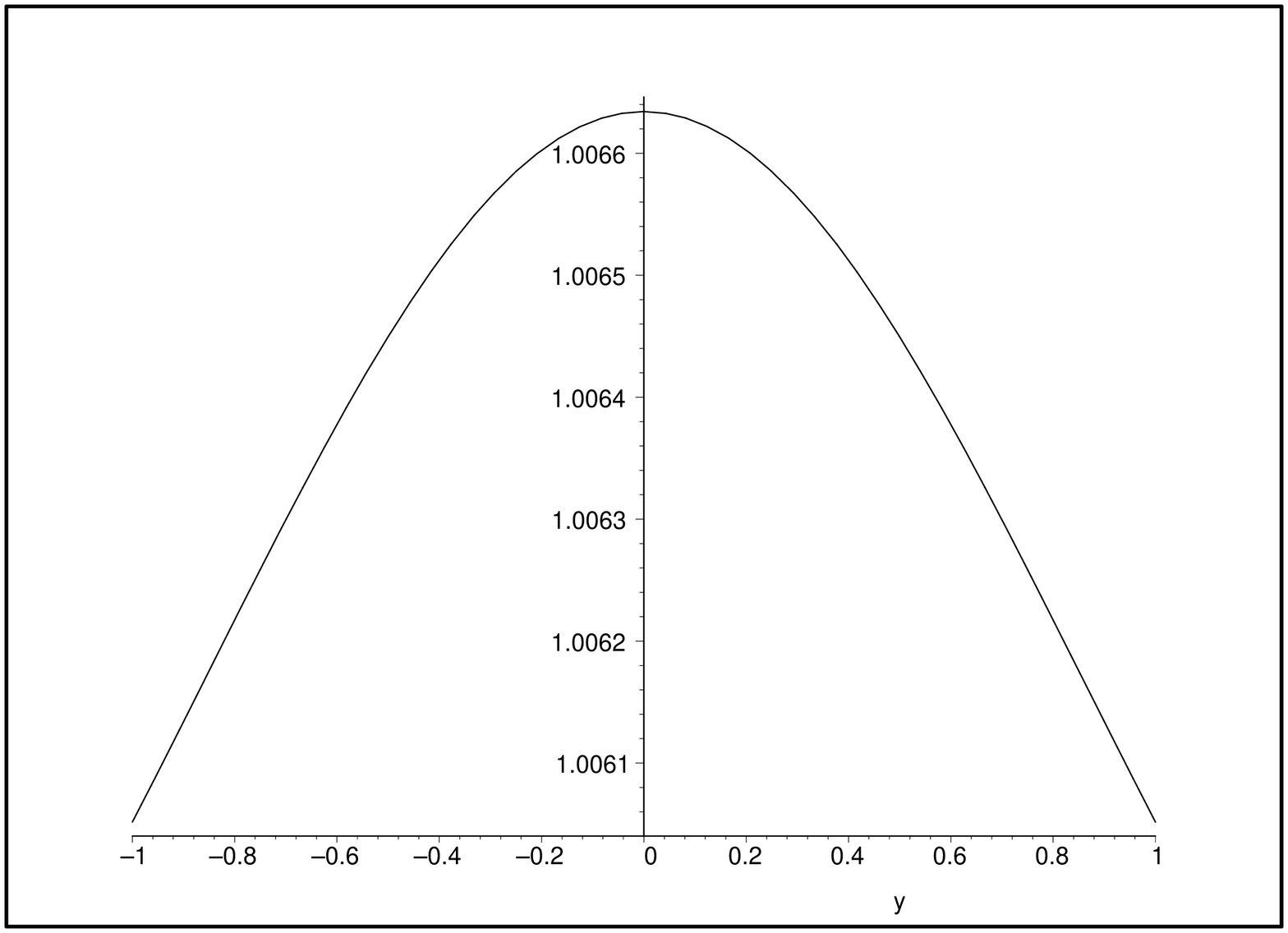}}}
\end{center}
\caption{A plot of the
ratio $\ln\left[  P_{40}\left(  \frac{y}{40}\right)  /40!\right]  /\ln\left[
\Psi_{2}(y,40)/40!\right]  $ .}%
\label{x=0}%
\end{figure}

We note from (\ref{Z1-}) that
\[
\sqrt{\frac{\pi}{2}}e^{x^{2}/2}w+Z_{1}(-x)=\sqrt{\frac{\pi}{2}}e^{x^{2}%
/2}\left(  w+2\right)  -Z_{1}(x),
\]
so that the integrand in (\ref{contour}) is antisymmetric with respect to the
map $\left(  x,w\right)  \rightarrow\left(  -x,-2-w\right)  .$ Thus, for
$n\rightarrow\infty$ and $x=O\left(  n^{-1}\right)  $ we get $P_{n}(x)\sim
\Psi_{2}(x,n)$ with%
\begin{equation}
\Psi_{2}(y,n)=\frac{n!}{\sqrt{2\ln\left(  n\right)  }}\left(  \sqrt{\frac
{2}{\pi}}\right)  ^{n+1}\left[  \exp\left(  y\sqrt{\frac{2}{\pi}}\right)
+\left(  -1\right)  ^{n}\exp\left(  -y\sqrt{\frac{2}{\pi}}\right)  \right]
.\label{18}%
\end{equation}
As $y\rightarrow\infty,$ the alternating term becomes negligible and
(\ref{18}) matches to (\ref{15}), as $x\rightarrow0.$ In Figure \ref{x=0} we plot the
ratio $\ln\left[  P_{40}\left(  \frac{y}{40}\right)  /40!\right]  /\ln\left[
\Psi_{2}(y,40)/40!\right]  $ and verify the accuracy of (\ref{18}).

Letting $x\rightarrow\infty$ in (\ref{14}) using (\ref{15}) yields%
\[
P_{n}(x)\sim\frac{n!\,x^{n+1}}{\sqrt{2\ln(n)}},
\]
which differs from (\ref{large x}). This suggests that another scale must be
analyzed, where $x$ and $n$ are both large. Thus, we consider the case of
$x\rightarrow\infty$, with $x=O\left(  \sqrt{\ln n}\right)  .$ Now the
singularity at $w=0$ in (\ref{contour}) becomes close to $w^{\ast}(x),$ since
\[
w^{\ast}(x)\sim-\frac{\sqrt{2}}{x\sqrt{\pi}},\quad x\rightarrow\infty.
\]
we use the form (\ref{Pn1}) and expand $S[N(x)+zN^{\prime}(x)]$ for
$z\rightarrow0$ and $x\rightarrow\infty.$ Setting $z=\xi/x$ with $\xi=O(1),$
we obtain
\begin{gather*}
S[N(x)+zN^{\prime}(x)]=\sqrt{2}\mathcal{J}\left\{  1+\sqrt{\frac{2}{\pi}%
}\left[  ze^{-x^{2}/2}-\zeta\left(  x\right)  \right]  \right\}  \\
=\sqrt{2}\mathcal{J}\left[  1+\sqrt{\frac{2}{\pi}}e^{-x^{2}/2}\left(
z-\frac{1}{x}+O\left(  x^{-3}\right)  \right)  \right]  \\
\sim\sqrt{2}\sqrt{-\ln\left[  \sqrt{\frac{2}{\pi}}e^{-x^{2}/2}\left(  \frac
{1}{x}-z\right)  \right]  }\\
\sim\sqrt{2}\sqrt{\frac{x^{2}}{2}+\ln(x)-\frac{1}{2}\ln\left(  \frac{2}{\pi
}\right)  -\ln\left(  1-\xi\right)  }.
\end{gather*}

Thus, we have%
\begin{equation}
P_{n}(x)\sim\frac{\sqrt{2}\left(  n+1\right)  !x^{n+1}}{2\pi\mathrm{i}}%
\oint\limits_{C_{1}}\sqrt{\frac{x^{2}}{2}+\ln\left(  \frac{x}{1-\xi}\right)
-\frac{1}{2}\ln\left(  \frac{2}{\pi}\right)  }\frac{d\xi}{\xi^{n+2}}.
\label{19}%
\end{equation}
Here the contour $C_{1}$ is a small loop about $\xi=0.$ Now we again employ
singularity analysis, with the branch point at $\xi=1$ determining the
asymptotic behavior for $n\rightarrow\infty.$ A deformation similar to that in
Figure \ref{Fig12} leads to%
\begin{equation}
P_{n}(x)\sim\frac{n!x^{n+1}}{\sqrt{2}}\frac{1}{\sqrt{\frac{x^{2}}{2}%
+\ln(nx)-\frac{1}{2}\ln\left(  \frac{2}{\pi}\right)  }}. \label{20}%
\end{equation}
For $x>>\sqrt{\ln(n)}$ this collapses to (\ref{large x}).

\begin{figure}[t]
\begin{center}
\rotatebox{270} {\resizebox{10cm}{!}{\includegraphics{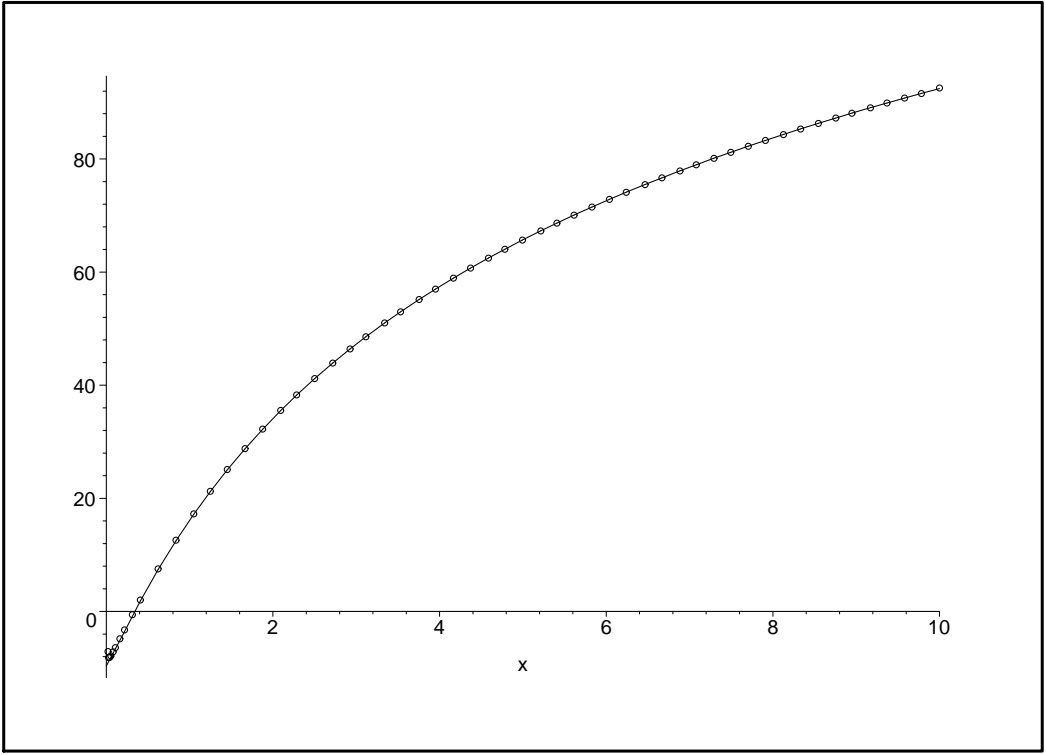}}}
\end{center}
\caption{A plot of $\ln\left[  P_{40}(x)/40!\right]  $ (solid line) and $\ln\left[  \Psi
_{3}(x,40)/40!\right]  $ (ooo).}%
\label{xlarge}%
\end{figure}

By examining (\ref{15}) and (\ref{20}), we can obtain the following
approximation%
\begin{equation}
P_{n}(x)\sim\Psi_{3}(x,n)=\frac{n!}{\sqrt{x^{2}+2\ln(nx)-\ln\left(  \frac
{2}{\pi}\right)  }}\left[  \frac{e^{-x^{2}/2}}{\zeta\left(  x\right)
}\right]  ^{n+1},\label{21}%
\end{equation}
which is more uniform in $x$, since it holds both for $x=O(1)$ and $x=O\left(
\sqrt{\ln n}\right)  $ for $n$ large and for $x\rightarrow\infty$ with $n$
fixed. However, we must still use (\ref{18}) if $n$ is large and $x$ is small.
In Figure \ref{xlarge} we plot $\ln\left[  P_{40}(x)/40!\right]  $ and $\ln\left[  \Psi
_{3}(x,40)/40!\right]  $ and confirm that (\ref{21}) is a better approximation
than (\ref{15}) for large values of $x.$ 

\section{WKB analysis\label{section3}}

We shall now rederive the results in the previous section by using only the
recurrence (\ref{2.3}) and (\ref{large x}). We apply the WKB method to
(\ref{2.3}), seeking solutions of the form $P_{n}(x)=n!\overline{P}_{n}(x),$
with%
\begin{equation}
\overline{P}_{n}(x)\sim\exp\left[  \left(  n+1\right)  A(x)\right]
B(x,n),\quad n\rightarrow\infty. \label{3.1}%
\end{equation}
Thus, we are assuming an exponential dependence on $n$ and an additional
weaker (e.g., algebraic) dependence that arises from the function $B(x,n).$
Using (\ref{3.1}) in (\ref{2.3}) leads to
\begin{align*}
&  e^{A(x)}\left[  B(x,n)+\frac{\partial}{\partial n}B(x,n)\right]  +O\left(
\frac{\partial^{2}B}{\partial n^{2}}\right) \\
&  =\left[  x+A^{\prime}(x)\right]  B(x,n)+\frac{1}{n+1}\frac{\partial
}{\partial x}B(x,n).
\end{align*}
Expecting that $\frac{\partial B}{\partial n}=o\left(  B\right)  $ and
$\frac{\partial^{2}B}{\partial n^{2}}=o\left(  \frac{\partial B}{\partial
n}\right)  ,$ we set%
\begin{equation}
e^{A(x)}=x+A^{\prime}(x) \label{3.3}%
\end{equation}
and%
\begin{equation}
e^{A(x)}\frac{\partial}{\partial n}B(x,n)=\frac{1}{n+1}\frac{\partial
}{\partial x}B(x,n)\sim\frac{1}{n}\frac{\partial}{\partial x}B(x,n).
\label{3.4}%
\end{equation}

To solve (\ref{3.3}) we let%
\[
A(x)=-\frac{x^{2}}{2}+a(x)
\]
to find that%
\[
a^{\prime}(x)=e^{-x^{2}/2}e^{a(x)}.
\]
Solving this separable ODE leads to%
\begin{equation}
A(x)=-\frac{x^{2}}{2}-\ln\left[  \zeta\left(  x\right)  +k\right]  ,
\label{3.7}%
\end{equation}
where $k$ is a constant of integration. To fix $k,$ we assume that expansion
(\ref{3.1}), as $x\rightarrow\infty,$ will asymptotically match to
(\ref{large x}), when this is expanded for $n\rightarrow\infty.$ In view of
(\ref{3.1}) this implies that%
\[
\overline{P}_{n}(x)\sim x^{n}=\exp\left[  n\ln(x)\right]  ,\quad
x\rightarrow\infty,
\]
so that%
\[
A(x)\sim\ln(x),\quad x\rightarrow\infty.
\]
In view of (\ref{3.7}) this is possible only if $k=0$ and then from
(\ref{zeta}) we have%
\begin{equation}
A(x)=-\ln\left[  e^{x^{2}/2}\zeta\left(  x\right)  \right]  \sim\ln(x),\quad
x\rightarrow\infty. \label{3.8}%
\end{equation}

We next analyze (\ref{3.4}). Using (\ref{3.8}) to compute $e^{A(x)},$ we
obtain%
\[
\frac{e^{-x^{2}/2}}{\zeta\left(  x\right)  }\frac{\partial}{\partial
n}B(x,n)=\frac{1}{n}\frac{\partial}{\partial x}B(x,n).
\]
Solving this first order PDE by the method of characteristics, we obtain%
\[
B(x,n)=b\left[  \frac{n}{\zeta\left(  x\right)  }\right]  ,
\]
where $b\left(  \cdot\right)  $ is at this point an arbitrary function.
However, since $n$ is large and $x=O(1),$ we need only the behavior of
$b\left(  \cdot\right)  $ for large values of its argument. We again argue
that by matching to (\ref{large x}) we have
\[
\exp\left[  \left(  n+1\right)  A(x)\right]  B(x,n)\sim x^{n},\quad
x\rightarrow\infty,
\]
and using (\ref{3.8}) we get%
\[
B(x,n)\sim e^{-A(x)}\sim\frac{1}{x},\quad x\rightarrow\infty
\]
and thus%
\[
b\left(  nxe^{x^{2}/2}\right)  \sim\frac{1}{x},\quad x\rightarrow\infty
\]
so that%
\[
b(z)\sim\frac{1}{\sqrt{2\ln(z)}},\quad z\rightarrow\infty.
\]

Combining our results, we have found that
\begin{equation}
P_{n}(x)\sim\frac{n!}{\sqrt{2}\sqrt{\ln(n)-\ln\left[  \zeta\left(  x\right)
\right]  }}\left[  \frac{e^{-x^{2}/2}}{\zeta\left(  x\right)  }\right]
^{n+1},\quad n\rightarrow\infty. \label{3.11}%
\end{equation}
This applies for $x=O(1)$ and $n\rightarrow\infty,$ where we can regain the
results of the singularity analysis (\ref{15}) by simply using%
\[
\sqrt{\ln(n)-\ln\left[  \zeta\left(  x\right)  \right]  }\sim\sqrt{\ln
(n)},\quad n\rightarrow\infty.
\]
Formula (\ref{3.11}) is also valid for $n=O(1)$ and $x\rightarrow\infty,$
where it reduces to (\ref{large x}). However, (\ref{3.1}) breaks down for
$x\rightarrow0.$

We thus consider the scale $x=y/n,$ $y=O\left(  1\right)  $ and set%
\begin{equation}
P_{n}(x)=n!\widetilde{P}_{n}(nx)=n!\widetilde{P}_{n}(y), \label{3.12}%
\end{equation}
with which (\ref{2.3}) becomes%
\begin{equation}
\widetilde{P}_{n+1}\left(  y+\frac{y}{n}\right)  =\frac{n}{n+1}\widetilde
{P}_{n}^{\prime}(y)+\frac{y}{n}\widetilde{P}_{n}(y). \label{3.13}%
\end{equation}
For fixed $y,$ we seek an asymptotic solution of (\ref{3.13}) in the form%
\begin{equation}
\widetilde{P}_{n}(y)\sim e^{\alpha n}q\left(  y,n\right)  ,\quad
n\rightarrow\infty, \label{3.14}%
\end{equation}
where $q\left(  y,n\right)  $ will have a weaker (e.g., algebraic or
logarithmic) dependence on $n.$ From (\ref{3.13}) we obtain, using
(\ref{3.14}),
\begin{align}
&  e^{\alpha}\left[  q(y,n)+\frac{\partial}{\partial n}q(y,n)+\frac{y}{n}%
\frac{\partial}{\partial y}q(y,n)+O\left(  n^{-2}\right)  \right]
\label{3.15}\\
&  =\left[  1-\frac{1}{n}+O\left(  n^{-2}\right)  \right]  \frac{\partial
}{\partial y}q(y,n)+\frac{y}{n}q(y,n).\nonumber
\end{align}
If $q\left(  y,n\right)  $ has an algebraic dependence on $n$, then
$\frac{\partial}{\partial n}q(y,n)$ should be roughly $O\left(  n^{-1}\right)
$ relative to $q\left(  y,n\right)  ,$ $\frac{\partial^{2}}{\partial n^{2}%
}q(y,n)$ roughly $O\left(  n^{-2}\right)  $ and so on. Thus, we expand
$q\left(  y,n\right)  $ as%
\begin{equation}
q\left(  y,n\right)  =q_{0}(y,n)+\frac{1}{n}q_{1}(y,n)+O\left(  n^{-2}\right)
, \label{3.16}%
\end{equation}
where $q_{0}(y,n),q_{1}(y,n)$ have a very weak (e.g., logarithmic) dependence
on $n$ and balance terms in (\ref{3.15}) of order $O(1)$ and $O\left(
n^{-1}\right)  $ to obtain%
\begin{equation}
e^{\alpha}q_{0}(y,n)=\frac{\partial}{\partial y}q_{0}(y,n) \label{3.17}%
\end{equation}
and%
\begin{align}
&  e^{\alpha}\left[  q_{1}(y,n)+\frac{\partial}{\partial n}q_{1}%
(y,n)+y\frac{\partial}{\partial y}q_{0}(y,n)\right] \label{3.18}\\
&  =\frac{\partial}{\partial y}q_{1}(y,n)-\frac{\partial}{\partial y}%
q_{0}(y,n)-yq_{0}(y,n).\nonumber
\end{align}

Solving (\ref{3.17}) yields
\begin{equation}
q_{0}(y,n)=\exp\left(  e^{\alpha}y\right)  \mathfrak{q}(n), \label{3.19}%
\end{equation}
where $\mathfrak{q}(n)$ must be determined. \ We could solve (\ref{3.18}),
using (\ref{3.19}), but its solution would involve another arbitrary function
of $n$. Thus, considering higher order terms will not help in determining
$\mathfrak{q}(n).$ Instead, we employ asymptotic matching to (\ref{3.11}).
Expanding (\ref{3.11}) for $x\rightarrow0$ and comparing the result to
(\ref{3.12}) as $y\rightarrow\infty,$ with (\ref{3.14}), (\ref{3.16}) and
(\ref{3.19}), we conclude that
\begin{equation}
\alpha=\frac{1}{2}\ln\left(  \frac{2}{\pi}\right)  ,\quad\mathfrak{q}%
(n)=\frac{1}{\sqrt{\pi\ln(n)}}. \label{3.20}%
\end{equation}

But then our approximation for $y=O(1)$ is not consistent with $P_{n}%
(0)=0=\widetilde{P}_{n}(0)$ for odd $n.$ We return to (\ref{3.13}) and observe
that the equation also admits an asymptotic solution of the form%
\[
\widetilde{P}_{n}(y)\sim\left(  -1\right)  ^{n}e^{\beta n}\overline{q}\left(
y,n\right)  ,\quad n\rightarrow\infty
\]
where, analogously to (\ref{3.17}), we find that
\[
-e^{\beta}\overline{q}_{0}(y,n)=\frac{\partial}{\partial y}\overline{q}%
_{0}(y,n),
\]
so that another asymptotic solution to (\ref{3.13}) is%
\begin{equation}
\widetilde{P}_{n}(y)\sim\left(  -1\right)  ^{n}e^{\beta n}\exp\left(
e^{\beta}y\right)  \overline{\mathfrak{q}}(n). \label{3.21}%
\end{equation}
We argue that any linear combination of (\ref{3.14}) and (\ref{3.21}) is also
a solution and that the combination which vanishes at $y=0$ for odd $n$ has
$\beta=\alpha$ and $\overline{\mathfrak{q}}(n)=\mathfrak{q}(n),$ as in
(\ref{3.20}). We have thus obtained, for $y=O(1),$
\[
P_{n}(x)\sim\frac{n!}{\sqrt{\pi\ln(n)}}\left(  \frac{2}{\pi}\right)
^{\frac{n}{2}}\left[  \exp\left(  y\sqrt{\frac{2}{\pi}}\right)  +\left(
-1\right)  ^{n}\exp\left(  -y\sqrt{\frac{2}{\pi}}\right)  \right]  ,\quad
n\rightarrow\infty.
\]
This agrees with (\ref{18}), obtained by singularity analysis in section 2.

To summarize, we have shown how to infer the asymptotics of $P_{n}(x)$ using
only the recursion (\ref{2.3}) and the large $x$ behavior (\ref{large x}). Our
analysis does need to make some assumptions about the forms of various
expansions and the asymptotic matching between different scales.

\section{The discrete ray method\label{section4}}

We shall now find a uniform asymptotic approximation for $P_{n}(x)$ using a
discrete form of the ray method \cite{MR0361328}. This approximation will
apply for $x$ and/or $n$ large. We seek an approximate solution for
(\ref{2.3}) of the form
\begin{equation}
P_{n}(x)=\exp\left[  f(x,n)+g(x,n)\right]  , \label{anszat}%
\end{equation}
where $g=o(f)$ as $n\rightarrow\infty.$ Since $P_{n}(x)=x^{n},$ $n=0,1$ we see
that we must have%
\begin{equation}
f(x,n)\sim n\ln(x)\text{ \ \ and \ \ }g(x,n)\rightarrow0 \label{initial}%
\end{equation}
as $n\rightarrow0.$ Using (\ref{anszat}) in (\ref{2.3}), we have%
\begin{equation}
\exp\left(  \frac{\partial f}{\partial n}+\frac{1}{2}\frac{\partial^{2}%
f}{\partial n^{2}}+\frac{\partial g}{\partial n}\right)  \sim\left(
\frac{\partial f}{\partial x}+\frac{\partial g}{\partial x}\right)  +\left(
n+1\right)  x \label{asymp1}%
\end{equation}
as $n\rightarrow\infty,$ where we have used%
\[
f(x,n+1)=f(x,n)+\frac{\partial f}{\partial n}(x,n)+\frac{1}{2}\frac
{\partial^{2}f}{\partial n^{2}}(x,n)+\cdots.
\]

From (\ref{asymp1}) we obtain, to leading order, the \textit{eikonal} equation%
\begin{equation}
\frac{\partial f}{\partial x}+\left(  n+1\right)  x-\exp\left(  \frac{\partial
f}{\partial n}\right)  =0, \label{eikonal}%
\end{equation}
and the \textit{transport }equation%
\begin{equation}
\frac{1}{2}\frac{\partial^{2}f}{\partial n^{2}}+\frac{\partial g}{\partial
n}-\frac{\partial g}{\partial x}\exp\left(  -\frac{\partial f}{\partial
n}\right)  =0. \label{transport}%
\end{equation}

To solve (\ref{eikonal}), we use the method of characteristics, which we
briefly review. Given the first order partial differential equation%
\[
F\left(  x,n,f,p,q\right)  =0,\text{ \ \ with \ \ }\ p=\frac{\partial
f}{\partial x},\quad q=\frac{\partial f}{\partial n},
\]
we search for a solution\ $f(x,n)$ by solving the system of \textquotedblleft
characteristic equations\textquotedblright\
\begin{align*}
\frac{dx}{dt}  &  =\frac{\partial F}{\partial p},\quad\frac{dn}{dt}%
=\frac{\partial F}{\partial q},\\
\frac{dp}{dt}  &  =-\frac{\partial F}{\partial x}-p\frac{\partial F}{\partial
f},\quad\frac{dq}{dt}=-\frac{\partial F}{\partial n}-q\frac{\partial
F}{\partial f},\\
\frac{df}{dt}  &  =p\frac{\partial F}{\partial p}+q\frac{\partial F}{\partial
q},
\end{align*}
with initial conditions%
\begin{equation}
F\left[  x(0,s),n(0,s),f(0,s),p(0,s),q(0,s)\right]  =0, \label{initial1}%
\end{equation}
and%
\begin{equation}
\quad\frac{d}{ds}f(0,s)=p(0,s)\frac{d}{ds}x(0,s)+q(0,s)\frac{d}{ds}n(0,s),
\label{initial2}%
\end{equation}
where we now consider $\left\{  x,n,f,p,q\right\}  $ to all be functions of
the variables $t$ and $s.$

For the eikonal equation (\ref{eikonal}), we have%
\begin{equation}
F\left(  x,n,f,p,q\right)  =p-e^{q}+\left(  n+1\right)  x\label{eikonal1}%
\end{equation}
and therefore the characteristic equations are%
\begin{equation}
\frac{dx}{dt}=1,\quad\frac{dn}{dt}=-e^{q},\quad\frac{dp}{dt}=-\left(
n+1\right)  ,\quad\frac{dq}{dt}=-x,\label{charac1}%
\end{equation}
and%
\begin{equation}
\frac{df}{dt}=p-qe^{q}.\label{eqf}%
\end{equation}
Solving (\ref{charac1}) subject to the initial conditions%
\[
x(0,s)=s,\quad n(0,s)=0,\quad p\left(  0,s\right)  =A(s),\quad q\left(
0,s\right)  =B(s),
\]
we obtain%
\begin{gather}
x=t+s,\quad n=-\sqrt{\frac{\pi}{2}}\exp\left(  \frac{s^{2}}{2}+B\right)
\left[  \operatorname{erf}\left(  \frac{t+s}{\sqrt{2}}\right)
-\operatorname{erf}\left(  \frac{s}{\sqrt{2}}\right)  \right]  ,\nonumber\\
p=\sqrt{\frac{\pi}{2}}\exp\left(  \frac{s^{2}}{2}+B\right)  \left(
t+s\right)  \left[  \operatorname{erf}\left(  \frac{t+s}{\sqrt{2}}\right)
-\operatorname{erf}\left(  \frac{s}{\sqrt{2}}\right)  \right]  \label{char0}\\
+\exp\left(  -\frac{1}{2}t^{2}-st+B\right)  -t-e^{B}+A,\quad\quad q=-\frac
{1}{2}t^{2}-st+B\nonumber
\end{gather}
From (\ref{initial}) we have%
\begin{equation}
A(s)=0\text{ \ and \ }B(s)=\ln(s),\label{A,B}%
\end{equation}
which is consistent with (\ref{initial1})$.$ Therefore,%
\begin{equation}
x=t+s,\quad n=-\sqrt{\frac{\pi}{2}}s\exp\left(  \frac{s^{2}}{2}\right)
\left[  \operatorname{erf}\left(  \frac{t+s}{\sqrt{2}}\right)
-\operatorname{erf}\left(  \frac{s}{\sqrt{2}}\right)  \right]  ,\label{rays}%
\end{equation}%
\begin{gather}
p=\sqrt{\frac{\pi}{2}}s\exp\left(  \frac{s^{2}}{2}\right)  \left(  t+s\right)
\left[  \operatorname{erf}\left(  \frac{t+s}{\sqrt{2}}\right)
-\operatorname{erf}\left(  \frac{s}{\sqrt{2}}\right)  \right]  \label{pq}\\
+s\exp\left(  -\frac{1}{2}t^{2}-st\right)  -\left(  t+s\right)  ,\quad\quad
q=-\frac{1}{2}t^{2}-st+\ln\left(  s\right)  .\nonumber
\end{gather}
In Figure \ref{rays} we sketch the rays $x(t,s),n(t,s)$ for $s\in\left[  -2..2\right]  .$

\begin{figure}[t]
\begin{center}
\rotatebox{270} {\resizebox{10cm}{!}{\includegraphics{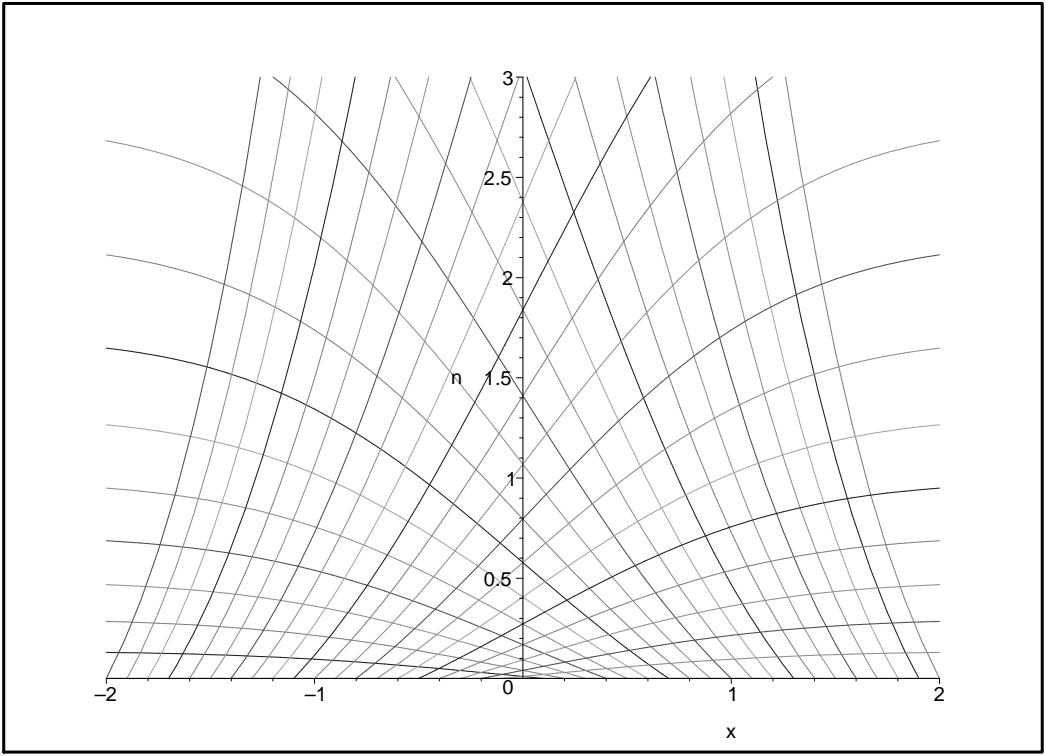}}}
\end{center}
\caption{A plot of the rays $x(t,s),n(t,s)$ for $s\in\left[  -2..2\right]  .$}%
\label{rays}%
\end{figure}

Using (\ref{pq}) in (\ref{eqf}) we have%
\begin{align}
\frac{df}{dt}  &  =\sqrt{\frac{\pi}{2}}s\exp\left(  \frac{s^{2}}{2}\right)
\left(  t+s\right)  \left[  \operatorname{erf}\left(  \frac{t+s}{\sqrt{2}%
}\right)  -\operatorname{erf}\left(  \frac{s}{\sqrt{2}}\right)  \right]
\label{eqf1}\\
&  +s\left[  1+\frac{1}{2}t^{2}+st-\ln\left(  s\right)  \right]  \exp\left(
-\frac{1}{2}t^{2}-st\right)  -\left(  t+s\right)  .\nonumber
\end{align}
Using (\ref{A,B}) in (\ref{initial2}), we get
\begin{equation}
f(0,s)=f_{0}, \label{f(0,s)}%
\end{equation}
and solving (\ref{eqf1}) subject to (\ref{f(0,s)}), we obtain%
\begin{align}
f(t,s)  &  =\sqrt{\frac{\pi}{2}}\exp\left(  \frac{s^{2}}{2}\right)  \left[
\operatorname{erf}\left(  \frac{t+s}{\sqrt{2}}\right)  -\operatorname{erf}%
\left(  \frac{s}{\sqrt{2}}\right)  \right] \label{f(t,s)}\\
&  \times s\left[  1+\frac{1}{2}t^{2}+st-\ln\left(  s\right)  \right]
-\left(  \frac{1}{2}t^{2}+st\right)  +f_{0}\nonumber
\end{align}
or, using (\ref{rays}),%
\begin{equation}
f=\left[  \ln\left(  s\right)  -1\right]  n+\frac{1}{2}\left(  s^{2}%
-x^{2}\right)  \left(  n+1\right)  +f_{0}. \label{f(x,n,s)}%
\end{equation}

To solve the transport equation (\ref{transport}), we need to compute
$\frac{\partial^{2}f}{\partial n^{2}},\frac{\partial g}{\partial n}$ and
$\frac{\partial g}{\partial x}$ as functions of $t$ and $s.$ Use of the chain
rule gives%
\[%
\begin{bmatrix}
\frac{\partial x}{\partial t} & \frac{\partial x}{\partial s}\\
\frac{\partial n}{\partial t} & \frac{\partial n}{\partial s}%
\end{bmatrix}%
\begin{bmatrix}
\frac{\partial t}{\partial x} & \frac{\partial t}{\partial n}\\
\frac{\partial s}{\partial x} & \frac{\partial s}{\partial n}%
\end{bmatrix}
=%
\begin{bmatrix}
1 & 0\\
0 & 1
\end{bmatrix}
\]
and hence,%
\begin{equation}%
\begin{bmatrix}
\frac{\partial t}{\partial x} & \frac{\partial t}{\partial n}\\
\frac{\partial s}{\partial x} & \frac{\partial s}{\partial n}%
\end{bmatrix}
=\frac{1}{J(t,s)}%
\begin{bmatrix}
\frac{\partial n}{\partial s} & -\frac{\partial x}{\partial s}\\
-\frac{\partial n}{\partial t} & \frac{\partial x}{\partial t}%
\end{bmatrix}
, \label{inversion}%
\end{equation}
where the Jacobian $J(t,s)$ is defined by
\begin{equation}
J\left(  t,s\right)  =\frac{\partial x}{\partial t}\frac{\partial n}{\partial
s}-\frac{\partial x}{\partial s}\frac{\partial n}{\partial t}=\frac{\partial
n}{\partial s}-\frac{\partial n}{\partial t}. \label{J}%
\end{equation}
Using (\ref{rays}), we can show after some algebra that%
\begin{equation}
J=\left(  s+\frac{1}{s}\right)  n+s. \label{J1}%
\end{equation}

Using $q=\frac{\partial f}{\partial n}$ in (\ref{transport}), we have%
\[
\frac{1}{2}\frac{\partial q}{\partial n}+\frac{\partial g}{\partial n}%
-\frac{\partial g}{\partial x}e^{-q}=0
\]
or%
\[
\frac{\partial}{\partial n}\left(  \frac{1}{2}e^{q}\right)  =\frac{\partial
g}{\partial x}-\frac{\partial g}{\partial n}e^{q}%
\]
and using (\ref{charac1}), we obtain%
\[
\frac{\partial}{\partial n}\left(  \frac{1}{2}e^{q}\right)  =\frac{\partial
g}{\partial x}\frac{\partial x}{\partial t}+\frac{\partial g}{\partial n}%
\frac{\partial n}{\partial t}=\frac{\partial g}{\partial t}.
\]
Since $-e^{q}=\frac{\partial n}{\partial t},$ we have%
\begin{gather*}
\frac{\partial}{\partial n}\left(  \frac{1}{2}e^{q}\right)  =-\frac{1}{2}%
\frac{\partial}{\partial n}\left(  \frac{\partial n}{\partial t}\right)
=-\frac{1}{2}\left(  \frac{\partial^{2}n}{\partial t^{2}}\frac{\partial
t}{\partial n}+\frac{\partial^{2}n}{\partial t\partial s}\frac{\partial
s}{\partial n}\right) \\
=-\frac{1}{2J}\left(  -\frac{\partial^{2}n}{\partial t^{2}}\frac{\partial
x}{\partial s}+\frac{\partial^{2}n}{\partial t\partial s}\frac{\partial
x}{\partial t}\right)  =-\frac{1}{2J}\left(  -\frac{\partial^{2}n}{\partial
t^{2}}+\frac{\partial^{2}n}{\partial t\partial s}\right) \\
=-\frac{1}{2J}\frac{\partial}{\partial t}\left(  \frac{\partial n}{\partial
s}-\frac{\partial n}{\partial t}\right)  =-\frac{1}{2J}\frac{\partial
J}{\partial t},
\end{gather*}
where we have used (\ref{inversion}) and (\ref{J}). Thus,%
\[
\frac{\partial g}{\partial t}=-\frac{1}{2J}\frac{\partial J}{\partial t}%
\]
and therefore%
\[
g(t,s)=-\frac{1}{2}\ln(J)+C(s)
\]
for some function $C(s).$ Since from (\ref{initial}) we have $g(0,s)=0,$ while
(\ref{J1}) gives $J(0,s)=s,$ we conclude that $C(s)=\frac{1}{2}\ln(s)$ and
hence
\begin{equation}
g(t,s)=\frac{1}{2}\ln\left[  \frac{s}{J(t,s)}\right]  . \label{g(t,s)}%
\end{equation}
Using (\ref{J1}) we can write (\ref{g(t,s)}) as%
\begin{equation}
g=\frac{1}{2}\ln\left[  \frac{s^{2}}{\left(  n+1\right)  s^{2}+n}\right]  .
\label{g(n,s)}%
\end{equation}
Replacing $f$ and $g$ in (\ref{anszat}) by (\ref{f(x,n,s)}) and (\ref{g(n,s)}%
), we obtain $P_{n}(x)\sim\Phi\left(  x,n;s\right)  $ as $n\rightarrow\infty,$
with
\begin{equation}
\Phi\left(  x,n;s\right)  =\kappa\exp\left[  \frac{1}{2}\left(  s^{2}%
-x^{2}\right)  \left(  n+1\right)  -n\right]  \frac{s^{n}}{\sqrt{n+1+ns^{-2}}%
}, \label{Phi}%
\end{equation}
where $\kappa=e^{f_{0}}$ is still to be determined.

Eliminating $t$ from (\ref{rays}) we get%
\begin{equation}
n+\sqrt{\frac{\pi}{2}}s\exp\left(  \frac{s^{2}}{2}\right)  \left[
\operatorname{erf}\left(  \frac{x}{\sqrt{2}}\right)  -\operatorname{erf}%
\left(  \frac{s}{\sqrt{2}}\right)  \right]  =0,\label{eqS}%
\end{equation}
which defines $s(x,n)$ implicitly. For every $n>0$ there exist only two
solutions $S_{m}(x,n)<0$ and $S_{p}(x,n)>0$ of (\ref{eqS}) (see Figure). Since
$\operatorname{erf}(x)$ is an odd function, it follows that
\begin{equation}
S_{m}(x,n)=-S_{p}(-x,n).\label{s(0,n)}%
\end{equation}
Although we have $P_{n}(x)\sim\Phi\left[  x,n;S_{m}(x,n)\right]  $ for
$x\ll-1$ and $P_{n}(x)\sim\Phi\left[  x,n;S_{p}(x,n)\right]  $ for $x\gg1,$
the two approximations are comparable when $x$ is small and therefore we must
add their contributions.

\begin{figure}[t]
\begin{center}
\rotatebox{270} {\resizebox{10cm}{!}{\includegraphics{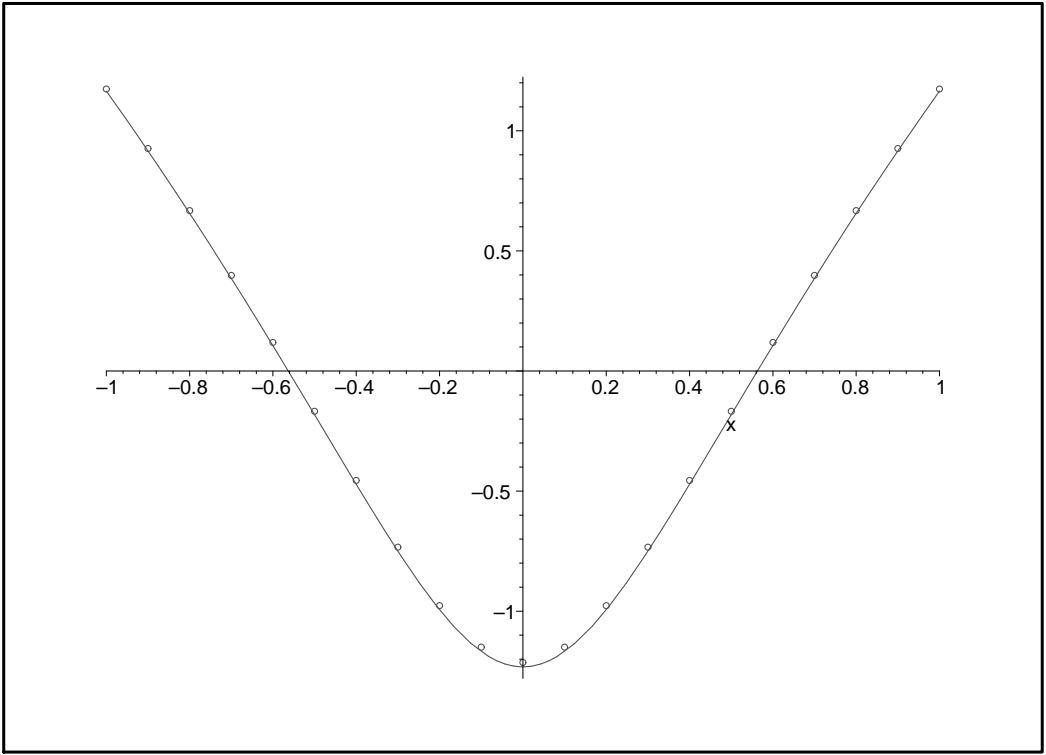}}}
\end{center}
\caption{A plot of $\ln\left[  P_{4}(x)/4!\right]  $ (solid line) and $\ln\left[
\Psi_{4}(x,4)/4!\right]  $ (ooo).}%
\label{ray0}%
\end{figure}

\begin{figure}[t]
\begin{center}
\rotatebox{270} {\resizebox{10cm}{!}{\includegraphics{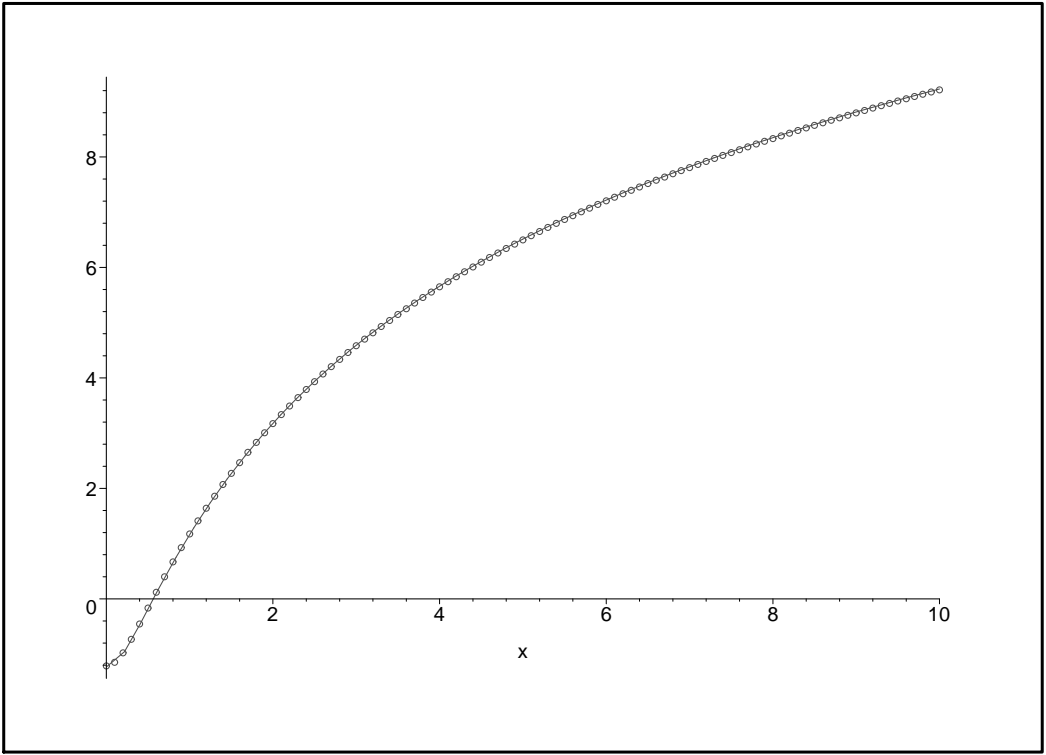}}}
\end{center}
\caption{A plot of $\ln\left[  P_{4}(x)/4!\right]  $ (solid line) and $\ln\left[
\Psi_{4}(x,4)/4!\right]  $ (ooo).}%
\label{raylarge}%
\end{figure}

We shall now find the constant $\kappa$ in (\ref{Phi}) by using (\ref{large x}%
). We rewrite (\ref{eqS}) as%
\begin{equation}
s^{2}\exp\left(  s^{2}\right)  =n^{2}\left[  \int\limits_{s}^{x}\exp\left(
-\frac{\theta^{2}}{2}\right)  d\theta\right]  ^{-2}\label{S1}%
\end{equation}
and for a fixed value of $n,$ consider the limit $x\rightarrow\infty.$ It
follows from (\ref{S1}) that $S_{p}(x,n)\sim x$ and therefore we consider an
expansion of the form%
\begin{equation}
S_{p}(x,n)\sim x+s_{0}+s_{1}x^{-1}+s_{2}x^{-2}+s_{3}x^{-3},\quad
x\rightarrow\infty.\label{S3}%
\end{equation}
Using (\ref{S3}) in (\ref{S1}), we obtain%
\begin{align*}
s_{0} &  =s_{2}=0,\quad s_{1}=\ln\left(  n+1\right)  ,\\
\quad s_{3} &  =1-\ln\left(  n+1\right)  -\frac{\ln^{2}\left(  n+1\right)
}{2}-\frac{1}{n+1}%
\end{align*}
and therefore%
\begin{equation}
s^{2}\exp\left(  s^{2}\right)  \sim\left(  n+1\right)  ^{2}x^{2}e^{x^{2}%
},\quad x\rightarrow\infty.\label{S4}%
\end{equation}
Solving (\ref{S4}) we have%
\begin{equation}
S_{p}(x,n)\sim\sqrt{\operatorname*{W}\left[  \left(  n+1\right)  ^{2}%
x^{2}e^{x^{2}}\right]  },\quad x\rightarrow\infty,\label{Spxlarge}%
\end{equation}
where $\operatorname*{W}\left(  z\right)  $ denotes the Lambert \textrm{W}
function, defined by \cite{MR1414285}%
\[
\operatorname*{W}\left(  z\right)  \exp\left[  \operatorname*{W}\left(
z\right)  \right]  =z,\quad\forall z\in\mathbb{C}%
\]
and having the asymptotic behavior%
\begin{equation}
\operatorname*{W}\left(  z\right)  =\ln(z)-\ln\ln\left(  z\right)  +\frac
{\ln\ln(z)}{\ln(z)}+O\left(  \left[  \frac{\ln\ln(z)}{\ln(z)}\right]
^{2}\right)  ,\quad z\rightarrow\infty.\label{Wasympt}%
\end{equation}
Using (\ref{Spxlarge}) and (\ref{Wasympt}) in (\ref{Phi}), we obtain%
\[
\Phi\left(  x,n;s\right)  \sim\kappa\left(  n+1\right)  ^{n+\frac{1}{2}}%
e^{-n}x^{n},\quad x\rightarrow\infty.
\]
From Stirling's formula%
\begin{equation}
n!=\left[  \sqrt{2\pi n}+O\left(  n^{-\frac{1}{2}}\right)  \right]
n^{n}e^{-n},\quad n\rightarrow\infty,\label{stirling}%
\end{equation}
and
\[
\left(  n+1\right)  ^{n+\frac{1}{2}}=\left[  e\sqrt{n}+O\left(  n^{-\frac
{1}{2}}\right)  \right]  n^{n},\quad n\rightarrow\infty,
\]
we conclude that%
\[
\kappa=e^{-1}\sqrt{2\pi}%
\]
and thus%
\begin{equation}
\Phi\left(  x,n;s\right)  =s^{n}\exp\left[  \frac{1}{2}\left(  s^{2}%
-x^{2}-2\right)  \left(  n+1\right)  \right]  \sqrt{\frac{2\pi s^{2}}{\left(
n+1\right)  s^{2}+n}}.\label{Phi1}%
\end{equation}
Using (\ref{s(0,n)}), we have $P_{n}(x)\sim\Psi_{4}(x,n)$ as $n\rightarrow
\infty,$ with%
\begin{equation}
\Psi_{4}(x,n)=\Phi\left[  x,n;S_{p}(x,n)\right]  +\Phi\left[  x,n;-S_{p}%
(-x,n)\right]  ,\quad n\rightarrow\infty.\label{Pfinal}%
\end{equation}
In Figure \ref{ray0} we compare $\ln\left[  P_{4}(x)/4!\right]  $ and $\ln\left[
\Psi_{4}(x,4)/4!\right]  $ for $0<x<10$ and in Figure \ref{raylarge} for $-1<x<1.$ We note
that the asymptotic approximation (\ref{Pfinal}) is more uniform than
(\ref{15}), (\ref{18}) and (\ref{20}) but it is less explicit since
$S_{p}(x,n)$ must be obtained numerically.

Next, we compare the results of this section with those in the previous two
sections. We first consider $x>0,$ with $x=O(1)$ and $n\rightarrow\infty.$
From (\ref{eqS}), we have%
\begin{equation}
S_{p}(x,n)\sim\sqrt{\operatorname*{W}\left[  \frac{\left(  n+1\right)  ^{2}%
}{\zeta^{2}(x)}\right]  }, \label{Snlarge}%
\end{equation}
where $\zeta(x)$ was defined in (\ref{zeta}). Using (\ref{Snlarge}) and
(\ref{Wasympt}) in (\ref{Phi1}), we get%
\[
P_{n}(x)\sim n^{n}e^{-n}\sqrt{\frac{n\pi}{\ln\left(  n\right)  }}\left[
\frac{e^{-x^{2}/2}}{\zeta(x)}\right]  ^{n+1},
\]
which agrees with (\ref{15}) after taking (\ref{stirling}) into account.

Now we consider the limit $n\rightarrow\infty,$ with $x=y/n$ and $y=O(1).$
From (\ref{eqS}), we have%
\begin{equation}
S_{p}\left(  y/n,n\right)  \sim\sqrt{\operatorname*{W}\left(  \frac{2n^{2}%
}{\pi}\right)  }+\left(  1+\sqrt{\frac{2}{\pi}}y\right)  \frac{1}%
{\sqrt{\operatorname*{W}\left(  \frac{2n^{2}}{\pi}\right)  }n}.\label{Sp}%
\end{equation}
Using (\ref{Sp}), (\ref{s(0,n)}) and (\ref{Wasympt}) in (\ref{Phi1}), we find
that%
\begin{align*}
\Phi\left(  y/n,n;S_{p}\right)   &  \sim n^{n}e^{-n}\sqrt{\frac{2n}{\ln(n)}%
}\left(  \sqrt{\frac{2}{\pi}}\right)  ^{n}\exp\left(  \sqrt{\frac{2}{\pi}%
}y\right)  ,\\
\Phi\left(  y/n,n;S_{m}\right)   &  \sim\left(  -1\right)  ^{n}n^{n}%
e^{-n}\sqrt{\frac{2n}{\ln(n)}}\left(  \sqrt{\frac{2}{\pi}}\right)  ^{n}%
\exp\left(  -\sqrt{\frac{2}{\pi}}y\right)
\end{align*}
and therefore%
\[
P_{n}(x)\sim n^{n}e^{-n}\sqrt{\frac{2n}{\ln(n)}}\left(  \sqrt{\frac{2}{\pi}%
}\right)  ^{n}\left[  \exp\left(  \sqrt{\frac{2}{\pi}}y\right)  +\left(
-1\right)  ^{n}\exp\left(  -\sqrt{\frac{2}{\pi}}y\right)  \right]
\]
agreeing with (\ref{18}).

Finally, we consider the limit $n\rightarrow\infty$ with $x=u\sqrt{\ln\left(
n\right)  },$ $u=O(1),$ $u>0.$ From (\ref{eqS}), we have%
\begin{equation}
S_{p}\left(  u\sqrt{\ln\left(  n\right)  },n\right)  \sim\sqrt
{\operatorname*{W}\left[  \left(  n+1\right)  ^{2}u^{2}n^{u^{2}}\ln(n)\right]
}. \label{S2}%
\end{equation}
Using (\ref{S2}) and (\ref{Wasympt}) in (\ref{Phi1}), we have%
\[
P_{n}\left(  u\sqrt{\ln\left(  n\right)  }\right)  \sim n^{n}e^{-n}\sqrt{2\pi
n}\frac{u^{n+1}}{\sqrt{u^{2}+2}}\left(  \sqrt{\ln\left(  n\right)  }\right)
^{n},\quad n\rightarrow\infty,
\]
which agrees with (\ref{20}).

\begin{acknowledgement}
The work of D. Dominici was supported by a Humboldt Research Fellowship for
Experienced Researchers from the Alexander von Humboldt Foundation. 

The work of C. Knessl was supported by the grants NSF 05-03745 and NSA H 98230-08-1-0102.
\end{acknowledgement}

\vspace*{0in}

\bibliographystyle{abbrv}

\begin{thebibliography}{10}

\bibitem{MR1225604}
M.~Abramowitz and I.~A. Stegun, editors.
\newblock {\em Handbook of mathematical functions with formulas, graphs, and
  mathematical tables}.
\newblock Dover Publications Inc., New York, 1992.

\bibitem{MR0153878}
L.~Carlitz.
\newblock The inverse of the error function.
\newblock {\em Pacific J. Math.}, 13:459--470, 1963.

\bibitem{MR1414285}
R.~M. Corless, G.~H. Gonnet, D.~E.~G. Hare, D.~J. Jeffrey, and D.~E. Knuth.
\newblock On the {L}ambert {$W$} function.
\newblock {\em Adv. Comput. Math.}, 5(4):329--359, 1996.

\bibitem{MR2364955}
D.~Dominici.
\newblock Asymptotic analysis of the {H}ermite polynomials from their
  differential-difference equation.
\newblock {\em J. Difference Equ. Appl.}, 13(12):1115--1128, 2007.

\bibitem{MR2401156}
D.~Dominici.
\newblock Asymptotic analysis of generalized {H}ermite polynomials.
\newblock {\em Analysis (Munich)}, 28(2):239--261, 2008.

\bibitem{MR2390273}
D.~Dominici.
\newblock Asymptotic analysis of the {K}rawtchouk polynomials by the {WKB}
  method.
\newblock {\em Ramanujan J.}, 15(3):303--338, 2008.

\bibitem{MR2250392}
D.~Dominici.
\newblock Asymptotic analysis of the {B}ell polynomials by the ray method.
\newblock {\em J. Comput. Appl. Math.}, (To appear.), 2009.

\bibitem{MR2117327}
D.~Dominici and C.~Knessl.
\newblock Geometrical optics approach to {M}arkov-modulated fluid models.
\newblock {\em Stud. Appl. Math.}, 114(1):45--93, 2005.

\bibitem{MR2262956}
D.~Dominici and C.~Knessl.
\newblock Ray solution of a singularly perturbed elliptic {PDE} with
  applications to communications networks.
\newblock {\em SIAM J. Appl. Math.}, 66(6):1871--1894 (electronic), 2006.

\bibitem{MR1986919}
D.~E. Dominici.
\newblock The inverse of the cumulative standard normal probability function.
\newblock {\em Integral Transforms Spec. Funct.}, 14(4):281--292, 2003.

\bibitem{MR0361328}
T.~Dosdale, G.~Duggan, and G.~J. Morgan.
\newblock Asymptotic solutions to differential-difference equations.
\newblock {\em J. Phys. A}, 7:1017--1026, 1974.

\bibitem{MR0434172}
V.~V. Frolov.
\newblock Group properties of the nonlinear heat-conduction equation and the
  solution of inverse problems.
\newblock {\em In\v z.-Fiz. \v Z.}, 30(3):546--553, 1976.

\bibitem{MR499726}
J.~B. Keller.
\newblock Rays, waves and asymptotics.
\newblock {\em Bull. Amer. Math. Soc.}, 84(5):727--750, 1978.

\bibitem{MR2078410}
C.~Knessl.
\newblock Some asymptotic results for the {$M/M/\infty$} queue with ranked
  servers.
\newblock {\em Queueing Syst.}, 47(3):201--250, 2004.

\bibitem{MR1695194}
C.~Knessl and W.~Szpankowski.
\newblock Quicksort algorithm again revisited.
\newblock {\em Discrete Math. Theor. Comput. Sci.}, 3(2):43--64 (electronic),
  1999.

\end{thebibliography}

\end{document}